\documentclass[a4paper, final, notitlepage, 10pt]{article}
\usepackage{latexsym,bm}
\usepackage{amsthm}
\usepackage{amsfonts}
\usepackage{amsmath}
\usepackage[all]{xy}
\usepackage[pdftex,final]{graphics}
\usepackage[dvipsnames]{pstricks}
\usepackage{fancyhdr}
\usepackage{verbatim}
\usepackage{geometry}
\usepackage{xcolor}
\usepackage{enumerate}
\usepackage{showkeys}
\usepackage{indentfirst}

\newtheorem{thm}{Theorem}[section]

\newtheorem{lem}[thm]{Lemma}
\newtheorem{prop}[thm]{Proposition}
\theoremstyle{definition}
\newtheorem{defin}[thm]{Definition}
\newtheorem{rem}[thm]{Remark}

\newtheorem{conj}[thm]{Conjecture}

\numberwithin{equation}{section}

\newcommand{\thmref}[1]{Theorem~\ref{#1}}

\newcommand{\lemref}[1]{Lemma~\ref{#1}}
\newcommand{\propref}[1]{Proposition~\ref{#1}}

\renewcommand{\O}{\mathcal{O}}
\renewcommand{\o}{\omega}

\title{Inequality of Noether type for Gorenstein minimal $3$-folds of general type}
\author{Yong Hu}
\date{}
\begin{document}
\maketitle
\renewcommand{\thefootnote}{\fnsymbol{footnote}}

\begin{abstract}
Let $X$ be a  Gorenstein minimal  $3$-fold of general type. We prove the optimal inequality:
\begin{align*}
K_X^{3}\geq \frac{4}{3}\chi(\omega_X)-2,
\end{align*}
where $\chi(\o_X)$ is the Euler-Poincar$\acute{\text{e}}$ characteristic of the dualizing sheaf $\o_{X}$.
\end{abstract}

\section{Introduction}
Throughout this paper, we work over the complex number field $\mathbb{C}$.

Let $S$ be a smooth minimal surface of general type. We have the classical Noether inequality: $K_S^2\ge 2p_g(S)-4$ and  $K_S^{2}\geq 2\chi(\O_S)-6$ (c.f.~\cite{N}). Together with the Miyaoka-Yau inequality $K_S^{2}\leq 9\chi(O_S)$ (cf.~\cite{Miyaoka}, \cite{Yau}), they have ever played very important roles in the classification of algebraic surfaces.

Let $X$ be a projective $3$-fold of general type. A natural question is: does there exists an inequality of Noether type for $3$-folds of general type?
There have been many works dedicated to proving the $3$-dimensional version of the Noether inequality:
\begin{enumerate}[\upshape$\cdot$]
\item In $1992$, M.~Kobayashi (c.f.~\cite[Proposition 3.2]{kobayashi}) constructed an infinite number of canonically polarized smooth $3$-folds of general type satisfying the equalities:
    \begin{align}\label{eq:examples}
    K_X^3=\frac{4}{3}p_g(X)-\frac{10}{3},\ K_X^3=\frac{4}{3}\chi(\omega_X)-2.
    \end{align}
\item In $2004$, M.~Chen (c.f.~\cite{Ch1}) studied minimal $3$-folds of general type and gave effective Noether type inequalities.
\item In $2004$, M.~Chen (c.f.~\cite{Ch2}) proved that the optimal inequality $K_X^3\ge \frac{4}{3}p_g(X)-\frac{10}{3}$ holds for all canonically polarized smooth $3$-folds of general type.
\item In $2006$, F.~Catanese, M.~Chen and De-Qi Zhang (c.f.~\cite{CCZ}) proved that the optimal inequality $K_X^3\ge \frac{4}{3}p_g(X)-\frac{10}{3}$ holds for all smooth minimal $3$-folds of general type.
\item In $2015$, J. A.~Chen and M.~Chen (c.f.~\cite{CC}) proved the optimal inequality $K_X^3\ge \frac{4}{3}p_g(X)-\frac{10}{3}$ under the assumption that $X$ is Gorenstein minimal.
\end{enumerate}

In this paper, a normal projective $3$-fold $X$ is called Gorenstein minimal if
 $X$ has at most $\mathbb{Q}$-factorial terminal singularities,
 the canonical divisor $K_X$ is a Cartier divisor and $K_X$ is nef.

 It is of interest to know whether there exists a similar Noether type inequality between $K_{X}^{3}$ and $\chi(\o_X)$. The following open problem was raised by M.~Chen (c.f.~\cite[3.9]{Ch2}):
\begin{conj}\cite[3.9]{Ch2}\label{conj:main conjecture}
Let $X$ be a Gorenstein minimal $3$-fold of general type. There should be an analogue of the Noether inequality in the form:
\begin{align*}
K_{X}^3\geq a\chi(\o_X)-b,
\end{align*}
where $a$ and $b$ are positive rational numbers.
\end{conj}
As was pointed out by M.~Chen (c.f.~\cite{Ch2}), it is difficult to find a Noether inequality in this direction because the inter relations among $p_g(X)$, $q(X)$ and $h^{2}(\mathcal{O}_X)$ are not clear to us, unlike in surface case. Some partial results were proved in \cite{Shin} and in \cite{CH1}.

\begin{enumerate}[\upshape$\cdot$]
\item In $1997$, D.~K.~Shin (c.f.~\cite{Shin}) proved that an effective inequality $K_X^3\ge\frac{6}{7}\chi(\o_X)-\frac{6}{7}$ holds for all smooth minimal $3$-folds of general type.
\item In $2006$, M.~Chen and C.~D.~Hacon (c.f.~\cite{CH1}) proved that an effective inequality $K_X^3\ge \frac{8}{9}\chi(\o_X)-\frac{10}{3}$ holds for all  smooth minimal $3$-folds of general type.
\end{enumerate}

We restrict our attention to the situation where $X$ is a  Gorenstein minimal $3$-fold of general type.
The aim of this paper is to prove the following.
\begin{thm}\label{thm:Main theorem}
Let $X$ be a  Gorenstein minimal $3$-fold of general type. Then
\begin{align}\label{eq:Noether inequality}
K_X^{3}\geq \frac{4}{3}\chi(\omega_X)-2.
\end{align}
\end{thm}
\begin{rem}
The inequality in \thmref{thm:Main theorem} is optimal because of M.~Kobayashi's examples (c.f.~\eqref{eq:examples}). It is well known that we have $\chi(\o_X)>0$ if $X$ is a Gorenstein minimal $3$-fold of general type. So  \eqref{eq:Noether inequality} is meaningful. One may ask whether \eqref{eq:Noether inequality} is still true if $X$ is a minimal $3$-fold of general type. Unfortunately, if $X$ is not Gorenstein, then $\chi(\o_X)$ could be either positive, zero, or negative (c.f.~\cite[line 10-18, page 2501]{CH2}). It is still make sense to ask whether a similar inequality holds for all minimal $3$-folds of general type. But this problem does not seem possible to resolve with the methods and the techniques of the present article.
\end{rem}
\section{Notations and the set up}

\begin{defin}\label{def:type}
Let $S$ be a smooth projective surface of general type. Denote by $S_0$ its minimal model and by $(a,b)=(K_{S_0}^2,p_g(S))$. We call $S$ is a surface of type $(a,b)$.
\end{defin}

Let $X$ be a Gorenstein minimal $3$-fold of general type with $p_g(X)\geq 2$. According to \cite[Lemma 5.1]{Kaw}, $X$ is locally factorial.
Write
\begin{align}\label{eq:movablefixed}
|K_X|=|\overline{M}|+\overline{Z},
\end{align}
where $|\overline{M}|$ is the movable part of $|K_X|$ and $\overline{Z}$ is the fixed part of $|K_X|$.

We shall resolve the base locus of $|\overline{M}|$ in two steps.
For a linear system $\Upsilon$, we denote by $\mathrm{Bs}\Upsilon$ the base locus of $\Upsilon$.
Roughly speaking,
the first step is to resolve the subset $\mathrm{Bs}|\overline{M}| \cap \mathrm{Sing}(X)$.
\begin{lem}[{cf.~\cite[Section~2]{CC}}]\label{lem:fgresolution}
There is a birational morphism $\alpha \colon X_0 \rightarrow X$
satisfying the following properties.
\begin{enumerate}[\upshape (a)]
\item The morphism $\alpha$ is a composition of successive divisorial contractions to points and $X_0$ is a Gorenstein $3$-fold with locally factorial terminal singularities.
\item Denote by $|M_0|$  the movable part of $|\alpha^*\overline{M}|$.
      Then $\mathrm{Bs}|M_0|\cap \mathrm{Sing}(X_0)=\emptyset$.
\item The following formulae
      \begin{align}
      K_{X_0}=\alpha^*K_X+\sum_{t=1}^{m}c_tD_t,\ \alpha^*(\overline{M})=M_0+\sum_{t=1}^{m}d_tD_t,\ \alpha^*(\overline{Z})=Z_0+\sum_{t=1}^{m}e_tD_t \label{eq:fGadjunction}
      \end{align}
      hold, where
      \begin{enumerate}[\upshape (i)]
      \item $Z_0$ is the strict transform of $\overline{Z}$,
         \item $D_t$ is a prime divisor such that $\alpha(D_t)$ is a point for $1 \le t \le m$, and
         \item $c_t$, $d_t$ and $e_t$ are non-negative integers such that $0 <c_t \le d_t$ for $1\le t \le m$.
      \end{enumerate}
\end{enumerate}
\end{lem}
\proof The birational morphism $\alpha$ is  constructed in \cite[p.~4--p.~5]{CC},
using explicit resolutions of terminal singularities (see \cite{jachen} and \cite[Definition~2.2]{CC}).
Then (a) and (b) follow from the construction and \cite[Lemma 5.1]{Kaw}.
Since both $X_0$ and $X$ are locally factorial, $c_t$, $d_t$ and $e_t$ are non-negative integers.
The inequality $c_t\le d_t$ follows by \cite[Corollary~2.4]{CC}.

\qed

We fix a birational morphism $\alpha \colon X_0 \rightarrow X$ as in \lemref{lem:fgresolution}.
We may assume that the number of divisorial contractions in the construction of $\alpha$ is minimal.
The second step is to resolve the base locus of $|M_0|$ without changing the singularities of $X_0$.
This is possible by \lemref{lem:fgresolution}~(b) and by Hironaka's Theorem (cf.~\cite{Hironaka}).
\begin{lem}[{cf.~\cite[Lemma~4.2]{Ch1}}]\label{lem:blowupresolution}
There are successive blowups
\begin{align*}
 \beta \colon Y=X_{n+1}\stackrel{\pi_n}\rightarrow X_{n}\rightarrow \cdots \rightarrow X_{i+1} \stackrel{\pi_{i}}\rightarrow X_{i}\rightarrow \cdots \rightarrow X_1\stackrel{\pi_0}\rightarrow X_0
\end{align*}
such that $\pi_i$ is a blowup along a smooth irreducible center $W_i$,
$W_i$ is contained in the base locus of the movable part of $|(\pi_0\circ \pi_1 \circ \cdots \circ \pi_{i-1})^*M_0|$ and $W_i \cap \mathrm{Sing}(X_i)=\emptyset$.
Moreover, the morphism $\beta=\pi_n\circ \cdots \circ \pi_0$ satisfies the following properties.
\begin{enumerate}[\upshape (a)]
\item Denote by $|M|$ the movable part of $|\beta^*M_0|$. Then $|M|$ is base point free.
\item The following formulae
      \begin{align}\label{eq:blowupadjunction}
          K_{Y}=\beta^*K_{X_0}+\sum_{i=0}^{n}a_i E_i,\ \beta^*M_0=M+\sum_{i=0}^{n}b_i E_i
      \end{align}
      hold, where $E_i$ is the strict transform of the exceptional divisor of $\pi_i$ for $0 \le i \le n$, and
      $a_i$ and $b_i$ are positive integers such that $a_i \le 2b_i\ \text{for}\ 0 \le i \le n$.
\end{enumerate}
\end{lem}
\proof
The construction of the blowups $\pi_i$ and  (a) follow by \lemref{lem:fgresolution}~(b) and by Hironaka's Theorem (cf.~\cite{Hironaka}).
We remark that the assertion $a_i \leq 2b_i$ in (b) is exactly  \cite[Lemma 4.2]{Ch1}.

\qed

From now on, we fix a birational morphism $\beta$ as in \lemref{lem:blowupresolution} such that the number $n+1$ of blowups is minimal.
Denote by $\phi_{K_X}$ the canonical map of $X$ and by $\Sigma$ the image of $\phi_{K_X}$.
Let $\phi$ be the morphism induced by the linear system $|M|$.
Then $\phi = \phi_{K_X} \circ \pi$, where $\pi=\alpha \circ \beta$.
Let $Y\stackrel{f}\rightarrow B \stackrel{\delta} \rightarrow \Sigma$ be the Stein factorization of $\phi$.
We have the following commutative diagram:
\begin{align*}
\xymatrix{
 & Y \ar[d]_{\pi} \ar[r]^{f} \ar[dr]^{\phi} \ar[dl]_{\beta}
                & B \ar[d]^{\delta} \\
X_0 \ar[r]^{\alpha} & X  \ar@{-->}[r]_{\phi_{K_X}}
                & \Sigma          }
\end{align*}
Note that $B$ is normal.
We have the following known results:
if $\dim B=3$, then $K_X^3\geq 2p_g(X)-6$ (cf.~\cite[Main Theorem]{kobayashi});
if $\dim B=2$, then $K_X^3\geq \left\lceil \frac{2}{3}(g(C)-1)\right\rceil(p_g(X)-2)$
where $g(C)$ is the genus of a general fiber $C$ of $f$ (cf.~\cite[Theorem 4.1~(ii)]{Ch1});
if $\dim B=1$, then either $K_X^3\geq 2p_g(X)-4$ or the general fiber of $f$ is a smooth projective surface of type $(1,2)$
(cf.~\cite[Theorem 4.1~(iii)]{Ch1}).

\section{Proof of \thmref{thm:Main theorem}}
This section is devoted to proving \thmref{thm:Main theorem}. Throughout this section, we denote by $X$ a Gorenstein minimal $3$-fold of general type. It is well known that $K_X^3$ is a positive even integer and $\chi(\omega)>0$(c.f.~\cite[2.1,\ 2.2]{JCCZ}). Denote by $a\colon X\rightarrow T$ the Stein factorization of the Albanese morphism of $X$ and by $F$ a general fiber of $a$. Since $3$-dimensional terminal singularities are isolated (c.f.~\cite{Reid}), $F$ is smooth.

The following lemma is due to \cite[Propsition 2.1]{CH1}.
\begin{lem}\label{lem:compare}
Let $X$ be a Gorenstein minimal $3$-fold of general type with $p_g(X)>0$. Then $\chi(\o_X)\le p_g(X)$ unless a general fiber of $a$ is a surface $F$ with $q(F)=0$, in which case one has the inequality
\begin{align*}
\chi(\o_X)\le (1+\frac{1}{p_g(F)})p_g(X).
\end{align*}
\end{lem}
\proof
After taking a resolution of $X$, \lemref{lem:compare} follows easily by \cite[Propsition 2.1]{CH1}.

\qed

We first prove \thmref{thm:Main theorem} under the assumption that $q(X)=0$.

\begin{lem}\label{lem:regular}
Let $X$ be a Gorenstein minimal $3$-fold of general type. Assume $q(X)=0$. Then
\begin{align}\label{eq:canonical volume reg}
K_X^3\geq \frac{4}{3}\chi(\o_X)-2.
\end{align}
\end{lem}
\proof
We have $\chi(\o_X)=p_g(X)-h^2(\mathcal{O}_X)-1\leq p_g(X)-1$ because $q(X)=0$.
Therefore \eqref{eq:canonical volume reg} holds by \cite[Theorem 1.1]{CC}.

\qed


 \begin{prop}\label{prop:dimY geq 2}
 Let $X$ be a Gorenstein minimal $3$-fold of general type. Keep the same notation as in the beginning of this section. Assume that $X$ satisfies one of the following conditions:
 \begin{enumerate}[\upshape (1)]
 \item $p_g(X)=0$.
 \item $q(X)=1$.
 \item $\mathrm{dim}T\ge 2$ and $p_g(X)>0$.
 \item $\mathrm{dim}T=1$ and $p_g(X)=1$.
 \end{enumerate}
 Then
  \begin{align}\label{eq:Noether inequalty in the first case}
  K_X^3\geq\frac{4}{3}\chi(\o_X)-\frac{10}{3}.
  \end{align}
 \end{prop}
 \proof
 If $p_g(X)=0$, we have $K_X^3\ge 2\chi(\o_X)> \frac{4}{3}\chi(\omega_X)-\frac{10}{3}$ by \cite[line 19-33, page 761]{CH1}. If $q(X)=1$, we have $\chi(\o_X)\leq p_g(X)$. Then \eqref{eq:Noether inequalty in the first case} holds by \cite[Theorem 1.1]{CC}. If $X$ satisfies condition (3), we have $\chi(\o_X)\le p_g(X)$ by \lemref{lem:compare}. So \eqref{eq:Noether inequalty in the first case} follows by \cite[Theorem 1.1]{CC}. Assume that $X$ satisfies condition (4). Since $p_g(X)>0$ and $F$ is a general fiber, we have $p_g(F)>0$. Then we can conclude $K_X^3\geq 2\chi(\o_X)> \frac{4}{3}\chi(\omega_X)-\frac{10}{3}$ from \lemref{lem:compare} and the fact that $K_X^3$ is a positive even integer.

 \qed

We now turn to the case where $X$ satisfies $p_g(X)\ge 2$, $q(X)\ge 2$ and $\mathrm{dim}T=1$. In this case, $T$ is a nonsingular projective curve with $g(T)=q(X)\ge 2$ (c.f.~\cite[Prop \uppercase\expandafter{\romannumeral 5}.15]{algebraicsurface}) and $F$ is a smooth minimal surface of general type because $X$ is minimal. The fibration $a$ is relatively minimal because $X$ is minimal. Therefore $K_{X/T}=K_X-a^*K_T$ is nef by \cite[Theorem 1.4]{O}.
Since $p_g(X)\ge 2$, we  can study the nontrivial canonical map of $X$. We can take the modification $\pi\colon Y\rightarrow X$ as in Section $2$ (see page $2$-$3$). Keep the same notation as in the last section. Recall that the  morphism $f\colon Y\rightarrow B$ is the Stein factorization of the canonical morphism (see page $3$).

\begin{lem}\label{lem:canoncial map not genus 2 fibration}
Let $X$ be a Gorenstein minimal $3$-fold of general type with $p_g(X)\geq 2$. Keep the same notation as above. Assume that $T$ is a nonsingular curve of genus $g(T)=q(X)\ge 2$.
Then
\begin{align}\label{eq:canonical map not genus 2 fib}
K_X^3\geq \frac{4}{3}\chi(\o_X)-6
\end{align}
unless the general fiber of $f$ is a curve of genus $g=2$.
\end{lem}
\proof
\textbf{Case $1$.} $\mathrm{dim}B\ge 2$.

From the discussion at the end of the last section, we have $K_X^{3}\geq 2p_g(X)-6$. Then we have \eqref{eq:canonical map not genus 2 fib} by \lemref{lem:compare}.

\textbf{Case $2$.} $\mathrm{dim}B=1$.

If the general fiber of $f$ is not a surface of type $(1,2)$, then we have $K_X^3\geq 2p_g(X)-4$ by the discussion at the end of the last section. So \eqref{eq:canonical map not genus 2 fib} holds by \lemref{lem:compare}.

Now we turn to the case where the general fiber of $f$ is a surface of type $(1,2)$. Since $K_X^3$ is a positive even integer, \eqref{eq:canonical map not genus 2 fib} holds if $\chi(\o_X)\leq 6$. So we may assume that $\chi(\o_X)\geq 7$. We have $p_g(X)\geq 5$ by \lemref{lem:compare}. We have $q(X)\le 1$ by \cite[Lemma 4.5]{Ch1}. But this contradicts to our assumption $q(X)\ge 2$. We are done.

\qed

\begin{prop}\label{prop:dim T=1 and g=2}
Let $X$ be a Gorenstein minimal $3$-fold of general type with $p_g(X)\geq 2$. Keep the same notation as above. Assume that $T$ is a nonsingular curve of genus $g(T)=q(X)\ge 2$ and that a general fiber $C$ of $f$ is a curve of genus $g(C)=2$. Then
\begin{align}\label{eq:Noether inequality fiber=2}
K_X^{3}\geq \frac{4}{3}\chi(\omega_X)-6.
\end{align}

\end{prop}
\proof
Since $K_X^3$ is a positive even integer, \eqref{eq:Noether inequality fiber=2} is automatically true for $\chi(\o_X)\leq 6$. We may assume  $\chi(\o_X)\geq 7$. So $p_g(X)\geq 5$ by \lemref{lem:compare}. According to \propref{prop:dimY geq 2}, we can assume that $q(X)\geq 2$ and that the general fiber $C$ of $f$ is a nonsingular curve of genus $2$. Recall that $|M|$ is the movable part of $|K_Y|$ and that $|M|$ is base point free.

We have $M^2 \equiv d_{\Sigma}\cdot \deg\delta \cdot C$, where $\Sigma$ is the image of $\phi$ and the symbol $\equiv$ stands for numerical equivalence.

Because $\Sigma$ is non-degenerate,  we have $d_{\Sigma} \ge p_g(X)-2$.
Since both $\pi^*K_X$ and $M$ are nef, we conclude that
\begin{align}\label{eq:canonical volume esti}
K_X^3 \ge (\pi^*K_X\cdot M^2) =d_{\Sigma}\cdot \deg\delta \cdot (\pi^*K_X\cdot C) \ge (p_g(X)-2)\deg\delta \cdot (\pi^*K_X\cdot C).
\end{align}
If $(\pi^*K_X\cdot C)\ge 2$,  then $K_X^3 \ge 2p_g(X)-4$. So \eqref{eq:Noether inequality fiber=2} holds by \eqref{eq:canonical volume esti} and \lemref{lem:compare}.

From now on, we  assume that $(\pi^*K_X\cdot C)=1$.

In order to prove \propref{prop:dim T=1 and g=2}, we need to use the techniques in the proofs of \cite[Theorem~4.3]{Ch1} and \cite[Theorem~3.1]{CC}. Recall from \eqref{eq:movablefixed}, \eqref{eq:fGadjunction} and \eqref{eq:blowupadjunction} that

\begin{align}\label{eq:fulladjunction}
K_Y=\pi^*K_X+(\sum_{t=1}^{m}c_t\beta^*D_t+\sum_{i=0}^{n}a_iE_i),&&
\pi^*K_X=M+(\sum_{t=1}^{m}(d_t+e_t)\beta^*D_t+\sum_{i=0}^{n}b_iE_i+\beta^*{Z_0})
\end{align}
In particular, we have
$K_X^3=(\pi^*K_X)^3=((\pi^*K_X)^2\cdot M)+(K_X^2\cdot\overline{Z})$.
We aim to bound $((\pi^*K_X)^2\cdot M)$ from below.

For this purpose, by Bertini's theorem, we choose a general member $S$ of $|M|$ such that $S$ is smooth and consider the fibration $f|_S$.
By abuse of notation, we still denote by $C$ the general fiber of $f|_S$.
We remark that $S$ is of general type since so is $X$.
Also the divisors $\beta^*D_t|_S$ and $E_i|_S$ are effective for $1 \le t \le m$ and $0 \le i \le n$.

Recall that $a\colon X\rightarrow T$ is the Stein factorization of the Albanese morphism of $X$ and that $K_{X/T}$ is a nef divisor. It is easy to see that $a_Y\colon Y\rightarrow T$ is the Stein factorization of the Albanese morphism of $Y$. Take a general fiber $F$ of $a$ such that $F_Y=\pi^*F$ is also a general fiber of $a_Y$. Note that $F$ is a smooth minimal surface of general type.

Since $\mathrm{dim}B=2$, we can see that $\mathrm{dim}f(F_Y)\geq 1$. So $S\cap F_Y$ is a nontrivial effective curve on $Y$. Thus we have $a_Y(S)=T$.

Since both $K_{X/T}=K_X-2(q(X)-1)F$ and $F_Y$ are nef divisors, we have $(\pi^*K_X\cdot C)\geq 2(q(X)-1)(F_Y\cdot C)\geq 0$.

Therefore $(F_Y\cdot C)=0$ because $q(X)\geq 2$ and $(\pi^*K_X\cdot C)=1$. So $f|_S$ and $a_Y|_S$ induce the same fibration. Denote by $\gamma\colon S\rightarrow \widehat{T}$ the induced fibration. We have $g(\widehat{T})\geq g(T)=q(X)\geq 2$.

By the definition of $M$ we can write
\begin{align}\label{eq:adjunction on S}
M|_S=\gamma^*A\equiv (\mathrm{deg}A)C,
\end{align}
where $A$ is an effective divisor on $\widehat{T}$.

Since $g(C)=2$,  we have $(K_Y\cdot C)=2$ by the adjunction formula.

According to \eqref{eq:fulladjunction} and $(\pi^*K_X\cdot C)=1$,
we have
\begin{align*}
((\sum_{t=1}^{m}c_t\beta^*D_t+\sum_{i=0}^{n}a_iE_i)|_S\cdot C)=1\ \text{and}\
((\sum_{t=1}^{m}(d_t+e_t)\beta^*D_t+\sum_{i=0}^{n}b_iE_i+\beta^*{Z_0})|_S\cdot C)=1.
\end{align*}
We conclude that the horizontal part of $(\sum_{t=1}^{m}c_t\beta^*D_t+\sum_{i=0}^{n}a_iE_i)|_S$ consists of an irreducible reduced curve $\Gamma$, which is also the horizontal part of $(\sum_{t=1}^{m}(d_t+e_t)\beta^*D_t+\sum_{i=0}^{n}b_iE_i+\beta^*{Z_0})|_S$ and $(\Gamma\cdot C)=1$.

We can write
\begin{align}
(\sum_{t=1}^{m}c_t\beta^*D_t+\sum_{i=0}^{n}a_iE_i)|_S&=\Gamma+D_V+E_V,\nonumber\\
(\sum_{t=1}^{m}(d_t+e_t)\beta^*D_t+\sum_{i=0}^{n}b_iE_i+\beta^*{Z_0})|_S&=\Gamma+D_V'+E_V'\label{eq:vertical}
\end{align}
where $E_V, D_V, E_V'$ and $D_V'$ are effective divisors contained in the fibers of $\gamma$.

According to \lemref{lem:fgresolution}~(c) and \lemref{lem:blowupresolution}~(b), we have
\begin{align}\label{eq:vertical comp}
D_V+E_V\leq 2D_V'+2E_V'.
\end{align}

Because $\Gamma$ is a section of $\gamma$,
\begin{align}\label{eq:gammaEV}
(\Gamma\cdot(2D_V'+2E_V'-D_V-E_V))\geq 0.
\end{align}

The adjunction formula yields
\begin{align}\label{eq:gammagenus}
(K_S\cdot\Gamma)+\Gamma^{2}=2p_a(\Gamma)-2 \ge 2q(X)-2
\end{align}
Note that
\begin{align*}
K_Y|_S=\pi^*K_X|_S+\Gamma+D_V+E_V\ \text{and}\ \pi^*K_X|_S=M|_S+\Gamma+D_V'+E_V'
\end{align*}
by \eqref{eq:fulladjunction}.
By \eqref{eq:gammaEV}, one has
\begin{align}
2q(X)-2 \le ((K_S+\Gamma)\cdot\Gamma)&=((K_Y+M)|_S\cdot\Gamma)+\Gamma^2 \nonumber\\
                  &=((\pi^*K_X|_S+M|_S+D_V+E_V+2\Gamma)\cdot\Gamma) \nonumber\\
                  &\le ((\pi^*K_X|_S+M|_S+2D_V'+2E_V'+2\Gamma)\cdot\Gamma) \label{eq:gammaEVeq}\\
                  &=((3\pi^*K_X|_S-M|_S)\cdot\Gamma) \nonumber\\
                  &=3(\pi^*K_X|_S\cdot\Gamma)-\mathrm{deg}A.\nonumber
\end{align}
The last equality holds by $M|_S \equiv  (\mathrm{deg}A) C$ (see \eqref{eq:adjunction on S}) and $(\Gamma\cdot C)=1$.

So we have
\begin{align}\label{eq:canonical degree of section}
(\pi^*K_X\cdot \Gamma)\geq \frac{1}{3}\mathrm{deg}A.
\end{align}
By \eqref{eq:adjunction on S} and $(\pi^*K_X\cdot C)=1$, we have
\begin{align}\label{eq:canoncial degree of vertical}
(\pi^*K_X|_S\cdot M|_S)=\mathrm{deg}A.
\end{align}

We have
\begin{align}\label{eq:canoncial volume ineq}
K_X^3\geq (\pi^*K_X|_S\cdot M|_S)+(\pi^*K_X\cdot\Gamma)\geq\frac{4}{3}\mathrm{deg}A
\end{align}
by \eqref{eq:canoncial degree of vertical} and \eqref{eq:canonical degree of section}.
We will bound $\mathrm{deg}A$ from below.

\textbf{Case $1$.} $h^1(\widehat{T}, A)>0$.

By Clifford's inequality and \eqref{eq:adjunction on S}, we have $$\mathrm{deg}A\geq 2h^0(\widehat{T},A)-2=2h^0(S,M|_S)-2\geq 2p_g(X)-4.$$
 So \eqref{eq:Noether inequality fiber=2} follows by above inequality, \eqref{eq:canoncial volume ineq} and \lemref{lem:compare}.

\textbf{Case $2$.} $h^1(\widehat{T}, A)=0$.

By Riemann-Roch formula and \lemref{lem:compare}, we have
\begin{align}\label{eq:degree of A}
\mathrm{deg}A=h^0(\widehat{T}, A)+g(\widehat{T})-1\geq p_g(X)+q(X)-2\geq\chi(\o_X)-1.
\end{align}
We conclude \eqref{eq:Noether inequality fiber=2} from \eqref{eq:canoncial volume ineq} and \eqref{eq:degree of A}.

\qed

We can now formulate our main result.

\begin{prop}\label{prop:irregular}
Let $X$ be an irregular Gorenstein minimal $3$-fold of general type. Then
\begin{align*}
K_X^3\geq\frac{4}{3}\chi(\o_X).
\end{align*}
\end{prop}
\proof
According to \propref{prop:dimY geq 2}, \lemref{lem:canoncial map not genus 2 fibration} and \propref{prop:dim T=1 and g=2}, we have
\begin{align}\label{eq:irregular}
K_X^3\geq \frac{4}{3}\chi(\o_X)-6.
\end{align}
Since $X$ is irregular, for every integer $m\geq2$ there is a cyclic unramified covering $\tau\colon \widehat{X}\rightarrow X$ of degree $m$. We have
 \begin{align}\label{eq:covering in the first case}
 K_{\widehat{X}}=\tau^*K_{X},\ K_{\widehat{X}}^3=mK_{X}^3,\ \chi(\o_{\widehat{X}})=m\chi(\o_X).
 \end{align}
 Note that $\widehat{X}$ is an irregular Gorenstein minimal $3$-fold of general type. So we have $K_{\widehat{X}}^3\geq \frac{4}{3}\chi(\o_{\widehat{X}})-6$. Therefore by \eqref{eq:irregular} and \eqref{eq:covering in the first case}, we have
 \begin{align*}
 K_X^3\geq\frac{4}{3}\chi(\o_X)-\frac{6}{m}.
 \end{align*}
 \propref{prop:irregular} follows by letting $m\rightarrow \infty$.

\qed

Thus  we can conclude \thmref{thm:Main theorem} from \lemref{lem:regular} and \propref{prop:irregular}.

\section{Acknowledgement}
The author would like to thank Professor Meng Chen for his guidance over this paper and encouragement. I wish to thank Yifan Chen for his helpful suggestions. The author is partially supported by the National Natural Science Foundation of China (Grant No.:~11571076).

\noindent\textbf{Authors' Addresses:}\\\smallskip

\noindent Yong~Hu,\\
School of Mathematical Sciences, Fudan University,\\
Shanghai 200433, P.~R.~China\\
Email:~yonghu11@fudan.edu.cn

\end{document}